\documentclass[10pt]{article}
\usepackage{amsmath,amssymb}
\usepackage{hyperref}
\usepackage{setspace}
\usepackage[utf8]{inputenc}
\usepackage[english]{babel}
\usepackage{amsfonts, amstext, amsthm, booktabs, dcolumn}
\usepackage{lscape}
\usepackage[english]{babel}
\usepackage[utf8]{inputenc}
\usepackage[T1]{fontenc}
\usepackage{authblk,color}
\usepackage[margin=1in]{geometry}

\newtheorem{theorem}{\bf Theorem}[section]

\newtheorem{corollary}{\bf Corollary}[section]

\newtheorem{remarks}{\bf Remarks}[section]
\newtheorem{lemma}{\bf Lemma}[section]

\newcommand\floor[1]{\left\lfloor#1\right\rfloor}

\begin{document}
\title{On Perturbations of Stein Operator}
\author[1]{A. N. Kumar}
\author[2]{N. S. Upadhye}
\affil[1]{\small Department of Mathematics, Indian Institute of Technology Madras,}
\affil[ ]{\small Chennai-600036, India.}
\affil[ ]{\small Email: amit.kumar2703@gmail.com}
\affil[2]{\small Department of Mathematics, Indian Institute of Technology Madras,}
\affil[ ]{\small Chennai-600036, India.}
\affil[ ]{\small Email: neelesh@iitm.ac.in}
\date{}
\maketitle

\begin{abstract}
\noindent
In this paper, we obtain Stein operator for sum of $n$ independent random variables (rvs) which is shown as perturbation of negative binomial (NB) operator. Comparing the operator with NB operator, we derive the error bounds for total variation distance by matching parameters. Also, three parameters approximation for such a sum is considered and is shown to improve the existing bounds in the literature. Finally, an application of our results to a function of waiting time for $(k_1,k_2)$-events is given. 
\end{abstract}
\begin{align*}
{\rm \bf Keywords:~}\text{Negative binomial distrib}&\text{ution, perturbation, probability generating function, Stein operator,} \\
&\hspace{-4cm}\text{ $(k_1,k_2)$ distribution.}\\
{\rm \bf MSC~2010~Subject~Classification:}&{~\rm Primary:~62E17,~62E20}\\
&\hspace{0.16cm}{\rm Secondary:~60E05,~60F05}
\end{align*}

\section{Introduction}
Applications of NB distribution appear in many areas such as network analysis, epidemics, telecommunications and related fields. NB approximation is widely studied in complex setting such as sum of waiting time, rare events and extremes. Also, NB approximation to sum of indicator rvs is given by Brown and Phillips \cite{BP}, approximation of NB and NB perturbation to the sum of the independent rvs is given by Vellaisamy et. al. \cite{VUC}, NB approximation to sum of independent NB rvs is given by Vellaisamy and Upadhye \cite{VN} and to $k$-runs is given by Wang and Xia \cite{WX}.\\ 
In this paper, we obtain a Stein operator for sum of independent rvs concentrated on ${\mathbb Z}_+=\{0,1,2,\dotsc\}$. Stein operator, so obtained, is perturbation of NB operator. So, we investigate error in approximation for NB to sum of independent rvs on ${\mathbb Z}_+$, via Stein method, by matching first and second moments. Also, error in approximation for convolution of NB and a geometric rv to sum of independent rvs on ${\mathbb Z}_+$ is investigated by matching first three moments. An application of these investigations is demonstrated for a function of waiting time for $(k_1,k_2)$-distribution. The approximation results, proved in Section \ref{AR} and Section \ref{app}, are either comparable to or improvement over the existing results in the literature. \\
Stein method (Stein \cite{stein}) is studied widely in probability approximations. For details and applications, see Barbour et. al. \cite{BCL}, Chen et. al. \cite{CGS}, Daly \cite{D1,D2}, Daly et. al. \cite{DLU}, Goldstein and Reinert \cite{GR}, Holmes \cite{SH}, Norudin and Peccati \cite{NP} and Ross \cite{R2}. For recent developments, see Barbour and Chen \cite{BAC}, Ley and Swan \cite{LS1,LS2}, Ley et. al. \cite{LRS}, Upadhye et. al. \cite{UCV} and references therein. This method involves identifying a suitable operator (known as a Stein operator) which can be obtained using one of the approaches (see Reinert \cite{R1}) such as, density approach (Stein \cite{stein,SDHR}), generator approach (Barbour and G\"{o}tze \cite{B,G}) and orthogonal polynomial approach (Diaconis and Zabell \cite{DZ}). Recently, probability generating function (PGF) approach (Upadhye et. al. \cite{UCV}) and a method to obtain canonical Stein operator (Ley et. al. \cite{LRS}) are developed. We focus on PGF approach for finding Stein operators.\\
The paper is organized as follows. In Section \ref{NAKR}, we define some necessary notations to formulate Stein method and our main results. Also, we explain some known results for NB distribution from the literature. In Section \ref{AR}, we first obtain Stein operator for sum of independent rvs which can be seen as perturbation of NB operator. So, we obtain bound between NB and sum of independent rvs by matching one and two parameters. Next, we derive Stein operator for convolution of NB and geometric which motivates us to use perturbation technique for obtaining bound between convolution of NB with a geometric and sum of independent rvs by matching three parameters. Finally, in Section \ref{app}, we give an application of our results for the function of waiting time for $(k_1,k_2)$-distribution.

\section{Notations and Known Results}\label{NAKR}
Throughout this paper, let $Z \sim {\rm NB}(\alpha, p)$ with 
$${\bf P}(Z = m) = {\alpha+m-1 \choose m} p^\alpha q^m, \quad \quad m = 0, 1, \dotsc,$$
where $\alpha > 0$ and $q = 1-p \in (0,1)$ and $Y = \sum_{i=1}^{n} X_i$, where $X_i,~i=1,2,\dotsc,n $ are independent rvs with PGF
\begin{equation}
M_Y(z) :={\mathbb E}\left(z^Y\right)= \sum_{m=0}^{\infty} {\bf P}(Y=m) z^m.\label{PGF}
\end{equation} 
Also, let the PGF of $X_i$ be $M_{X_i}$ such that 
\begin{equation}
\displaystyle{G_{X_i}(z) := \frac{M^\prime_{X_i}(z)}{M_{X_i}(z)} = \sum_{m=0}^{\infty} a_{i,m+1} z^m}.\label{assumption}
\end{equation} 
In particular, for specific distributions, following holds.
\begin{itemize}
\item[{\bf (O1)}] $X_i \sim Ge(p_i) \implies a_{i,m+1} = q_i^{m+1}.$
\item[{\bf (O2)}] $X_i \sim Bi(\tilde{n},\tilde{p}_i) \implies a_{i,m+1} = \tilde{n} (-1)^m \left(\tilde{p}_i/\tilde{q}_i\right)^{m+1}$.
\item[{\bf (O3)}] $X_i \sim Po(\lambda_i) \implies a_{i,m+1} = \lambda_i$ for $m=0$ and $0$ otherwise.
\end{itemize}
Next, let $\mu$ and $\sigma^2$ denote the mean and variance of $Y$ respectively. Then
\begin{equation}
\begin{aligned}
\mu &:= \sum_{i=1}^{n}G_{X_i}(1) = \sum_{i=1}^{n} \sum_{m=0}^{\infty} a_{i,m+1}, \quad \sigma^2 := \sum_{i=1}^{n}\left[G_{X_i}(1) + G^\prime_{X_i}(1)\right] = \sum_{i=1}^{n}\sum_{m=0}^{\infty} (m+1) a_{i,m+1},\\
\mu_2 &:= \sum_{i=1}^{n} G^\prime_{X_i}(1) = \sum_{i=1}^{n} \sum_{m=0}^{\infty} m a_{i,m+1} \quad {\rm and} \quad \mu_3 := \sum_{i=1}^{n}G^{\prime\prime}_{X_i}(1) = \sum_{i=1}^{n}\sum_{m=0}^{\infty} m(m-1) a_{i,m+1},\label{notation}
\end{aligned}
\end{equation}
where $\mu_2$ and $\mu_3$ denote second and third factorial cumulant moment of $Y$ (see Vellaisamy et. al. \cite{VUC}, p.p. $104-105$). Let us define
\begin{equation}
\begin{split}
\eta_1 := \frac{3}{2} \mu {\mu_3} - 4 {\mu_2}^2, \quad \eta_2 := 27 \mu^2 {\mu_2}^2-16 {\mu_2}^3 - \frac{27}{2} \mu^3 {\mu_3} + 9 \mu {\mu_2}{\mu_3},\\
\eta_3 := \left(\frac{\eta_2 + \sqrt{4 \eta_1^3 + \eta_2^2}}{2}\right)^{1/3} \quad {\rm and} \quad\eta := 2 {\mu_2} + \frac{\eta_1}{\eta_3} - \eta_3,\label{eta}
\end{split} 
\end{equation}
provided $\eta,\eta_3 \in {\mathbb R}_+$, the set of all positive real numbers. \\
Now, let $\cal G$ be the set of all bounded function on ${\mathbb Z}_+$ and 
\begin{equation}
{\cal G}_X = \left\{g|~g \in {\cal G}~{\rm such~that}~g(0)=0~{\rm and} ~g(x)=0, ~{\rm for}~x \notin {supp(X)}  \right\}\label{gx}
\end{equation} 
be associated with Stein operator ${\cal A}_X$, where $Supp(X)$ denotes the support of a rv $X$.\\
Next, Stein method can be formulated in three steps. First, identify a suitable operator (known as a Stein operator) for the rv $X$. Stein operator is defined on family of function ${\cal G}_X$ such that
$${\mathbb E}\left({\cal A}_X g\right) = 0, \quad \quad {\rm for}~g \in {\cal G}_X.$$
In second step, we find the solution (say $g_f$) of the difference equation (known as Stein equation)
\begin{equation}
{\cal A }_X g(m) = f(m) - {\mathbb E}f(Z), \quad m \in {\mathbb Z}_+~{\rm and}~f \in {\cal G}.\label{steq}
\end{equation}
and obtain the bound for $g_f$ (or $\Delta g_f$, as required) in terms of $f$.\\
Finally, Substituting a rv $Y$ for $m$ in (\ref{steq}) and taking expectations and supremum, we get the following
\begin{equation} 
d_{TV}\left(Y,X\right) := \sup_{f \in {\cal H}}|{\mathbb E}f(Y) - {\mathbb E}f(X)| = \sup_{f \in {\cal H}}|{\mathbb E}\left[{\cal A }_X g_f(Y)\right]|,\label{tv}
\end{equation}
where ${\cal H} = \{{\mathbb I}_S|~ S~measurable\}$ and ${\mathbb I}_S$ is the indicator function of the set $S$. Equation (\ref{tv}) is also equivalent to
$$d_{TV}(Y, X) = \frac{1}{2} \sum_{m=0}^{\infty} \left|{\bf P}(Y=m) - {\bf P}(X = m)\right|.$$
As $Y$ is sum of independent rvs on ${\mathbb Z}_+$, from Corollary $1.6$ of Mattner and Roos \cite{MR}, we have
\begin{equation}
d_{TV} (Y, Y+1) \le \sqrt{\frac{2}{\pi}} \left(\frac{1}{4} + \sum_{i=1}^{n}(1 - d_{TV}(X_i, X_{i}+1))\right)^{-1/2}.\label{mros}
\end{equation}
Next, it is known that Stein operator for NB$(\alpha, p)$ is given by (Brown and Phillips \cite{BP})
\begin{equation}
\left({\cal A}_Z g\right)(m) = q(\alpha + m)g(m+1) - m g(m), \quad \quad {\rm for}~m \in {\mathbb Z}_{+}~{\rm and}~g \in {\cal G}_Z. \label{NBoper}
\end{equation}
Also, bound for the solution to (\ref{steq}) is given by  
\begin{equation}
\|\Delta g\| \le 1/\alpha q,\label{bound}
\end{equation} 
where $\|\Delta g\|=\sup_{m \in {\mathbb Z}_+} |\Delta g(m)|$ and $\Delta g(m)= g(m+1) - g(m)$ denotes first forward difference operator (see Brown and Phillips \cite{BP} and Vellaisamy et. al. \cite{VUC} for details).\\
As NB distribution can be described using two parameters, namely $\alpha$ and $p$, we can study NB approximation problem using the moment matching technique up to first two moments. For an approximation with extra parameter, we can use perturbation technique described by Barbour et. al.  \cite{BCX} and can be formulated for NB distribution as follows:\\
Let ${\cal A}_V$ be a Stein operator of $V = Z + W$, where $W$ be a rv with parameter $\hat{p}=1-\hat{q}$ and probability mass function (PMF)
$${\bf P}(W=x) = \hat{q}^x \hat{p}, \quad \quad x=0,1,2, \dotsc.$$ 
Let ${\cal U}_V = {\cal A}_V - {\cal A}_Z$, such that
$$\|{\cal U}_V g\| \le \delta_1 \|\Delta g\|,$$
where $\alpha q > \delta_1$ and $g \in {\cal G}_V$, defined in (\ref{gx}). Also, let the rv $Y$ satisfy
$$\left|{\mathbb E}\left({\cal A}_V g\right)(Y)\right| \le \delta_2 \|\Delta g\|, \quad \quad  {\rm for} ~\delta_2 \ge 0$$
then
\begin{equation}
d_{TV}(Y,V) \le \frac{\delta_2}{\alpha q - \delta_1}. \label{perturbation}
\end{equation}
(See Theorem $2.4$ of Barbour et. al. \cite{BCX} and $(8)$ of Vellaisamy et. al. \cite{VUC} for more details).

\section{Approximation Results}\label{AR}
In this section, we obtain $Z$-approximation and $V$-approximation bounds to $Y$ using first two moments and first three moments respectively.

\subsection{One-parameter approximation}
The choice of the parameters can be done using the following relation.
\begin{equation}
\frac{\alpha q}{p} = \mu \quad \implies \quad p=\frac{\alpha}{\alpha+\mu}~~\left({\rm or}~\alpha = \frac{\mu p}{q}\right). \label{Mmean}
\end{equation}
Here, matching can be done in two ways:
\begin{itemize}
\item[(i)] Let $\alpha$ be fixed (in particular, $\alpha = n$) and $p=\alpha/(\alpha+\mu)$.
\item[(ii)] Let $p$ be fixed of our choice and the choice of $\alpha=\mu p/q$.
\end{itemize}

\begin{theorem}\label{fmth}
Let $X_1, X_2, \dotsc,X_n$ are independent rvs with (\ref{assumption}) and $Y=\sum_{i=1}^{n} X_i$, then 
$$d_{TV}\left(Y, Z\right) \le \frac{1}{\alpha q} \sum_{i=1}^{n} \sum_{l=1}^{\infty} l \left|a_{i,l+1} - q a_{i,l}\right|,$$
where $Z \sim {\rm NB}(\alpha,p)$.
\end{theorem}

\noindent
{\bf Proof}. Given $Y = \sum_{i=1}^{n} X_i$ such that $X_i,~i=1,2.\dotsc,n$ are independent rvs. Then the PGF of $Y$ is given by $M_{Y}(z) = \prod_{i=1}^{n} M_{X_i}(z)$. Differentiating with respect to $z$, we have 
$$M_Y^\prime(z) = M_Y(z) \sum_{i=1}^{n} G_{X_i}(z) = \sum_{i=1}^{n} M_Y(z) \left(\sum_{m=0}^{\infty}a_{i,m+1} z^m\right),$$
where $G_{X_i}(\cdot)$ as defined in (\ref{assumption}). Using (\ref{PGF}) and multiplying by $(1-q z)$, we get
$$\sum_{m=0}^{\infty}(m+1)p_{m+1}z^m - q \sum_{m=0}^{\infty}(m+1)p_{m+1}z^{m+1} = \sum_{i=1}^{n} \left[\sum_{m=0}^{\infty}\left(\sum_{l=0}^{m}p_l a_{i,m-l+1}\right)z^m - q \sum_{m=0}^{\infty}\left(\sum_{l=0}^{m-1} p_l a_{i,m-l}\right)z^{m}\right],$$
where $q = 1-p$ and $p$ is defined in (\ref{Mmean}). Now, comparing the coefficient of $z^m$, we obtain the recursive relation
$$q m p_m - (m+1)p_{m+1} + \sum_{i=1}^{n} \left(\sum_{l=0}^{m} p_l a_{i,m-l+1} - q \sum_{l=0}^{m-1}p_l a_{i,m-l}\right) = 0.$$
Let $g \in {\cal G}_Y$, defined in (\ref{gx}), then
$$\sum_{m=0}^{\infty} g(m+1)\left[q m p_m - (m+1)p_{m+1} + \sum_{i=1}^{n} \left(\sum_{l=0}^{m} p_l a_{i,m-l+1} - q \sum_{l=0}^{m-1}p_l a_{i,m-l}\right)\right] = 0,$$
or equivalently
$$\sum_{m=0}^{\infty}\left[q m g(m+1) - m g(m) + \left(\sum_{i=1}^{n} a_{i,1}\right) g(m+1) + \sum_{i=1}^{n}\sum_{l=1}^{\infty}g(l+m+1)\left(a_{i,l+1} - q a_{i,l}\right)\right]p_m = 0.$$
Hence, Stein operator for $Y$ is given by
$${\cal A}_Y g(m) =q m g(m+1) - m g(m) + \left(\sum_{i=1}^{n} a_{i,1}\right) g(m+1) + \sum_{i=1}^{n}\sum_{l=1}^{\infty}g(l+m+1)\left(a_{i,l+1} - q a_{i,l}\right).$$
Rewrite Stein operator, using auxiliary parameter $\alpha > 0$, as
\begin{equation}
{\cal A}_Y g(m) = q\left(\alpha+m\right)g(m+1) - m g(m) + \left(\sum_{i=1}^{n} a_{i,1} - \alpha q \right)g(m+1)+ \sum_{i=1}^{n}\sum_{l=1}^{\infty}g(l+m+1)\left(a_{i,l+1} - q a_{i,l}\right).\label{opertilde}
\end{equation}
This is a Stein operator for sum of independent rvs, which is a perturbation of NB$(\alpha, p)$ in view of Barbour and Xia \cite{BX} and Vellaisamy et. al. \cite{VUC}. Applying Newton's expansion as given in Barbour and \v{C}ekanavi\v{c}ius \cite{BC}, we have
\begin{equation}
g(m+l+1) = \sum_{j=1}^{l}\Delta g(m+j) + g(m+1). \label{delta}
\end{equation}
Putting (\ref{delta}) in (\ref{opertilde}) and  using (\ref{Mmean}), we get
\begin{align*}
{\cal A}_Y g(m) &= q(\alpha+m)g(m+1) - m g(m) + \sum_{i=1}^{n} \sum_{l=1}^{\infty} \sum_{j=1}^{l} \Delta g(m+j) \left(a_{i,l+1} - q a_{i,l}\right)&= {{\cal A}}_{Z} g(m) + {{\cal U}_Y} g(m),
\end{align*}
where ${{\cal A}}_{Z}$ is a Stein operator for NB ($\alpha,p$) described as in (\ref{NBoper}). ${\cal A}_Y$ is a Stein operator for sum of $n$ independent rvs by matching mean with negative binomial rv. Now, for $g \in {\cal G}_Z \cap {\cal G}_Y$, taking the expectation of ${{\cal U}_Y}$ with respect to $Y$ and using (\ref{bound}), we get required result.\qed

\begin{corollary}\label{cor1}
Given $Y=\sum_{i=1}^{n} X_i$, let $X_i$ are different type of distribution, we have the following bounds
\begin{itemize}
\item[{(i)}] Let $X_i$ follow Ge$(p_i),~i=1,2,\dotsc,n$ with $q_i = (1-p_i) < 1/2$, then
\begin{equation}
d_{TV}(Y,Z) \le \frac{1}{\alpha q} \sum_{i=1}^{n} \left|p-p_i\right| \sigma^2_{X_i},\label{gbop}
\end{equation} 
where $\sigma^2_{X_i}$ is the variance of $X_i$.
\item[{(ii)}] Let $X_i$ follow Po$(\lambda_i)$ for $i \in S_1$ and Ge$(p_i)$ for $i \in S_2$, where $S_1 \cup S_2 = \{1,2,\dotsc,n\}$, then
\begin{equation}
d_{TV}(Y,Z) \le \frac{1}{\alpha q} \left(q \sum_{i \in S_1} \lambda_i + \sum_{i \in S_2} |p-p_i| \frac{q_i}{p_i^2}\right).\label{pgop}
\end{equation}
\item[{(iii)}] Let $X_i$ follow Bi$(\tilde{n},\tilde{p}_i)$ for $i \in S_1$ and Ge$(p_i)$ for $i \in S_2$, where $S_1 \cup S_2 = \{1,2,\dotsc,n\}$ with $q_i, \tilde{p}_i < 1/2$, then
\begin{equation}
d_{TV}(Y,Z) \le \frac{1}{\alpha q} \left(\tilde{n} \sum_{i \in S_1} \left(\frac{\tilde{p}_i}{\tilde{q}_i}+q\right) \frac{\tilde{p}_i \tilde{q}_i}{(1-2\tilde{p}_i)^2}+q \sum_{i \in S_2} |p-p_i| \frac{q_i}{p_i^2} \right).\label{bgop}
\end{equation}
\item[{(iv)}] Let $X_i$ follow Po$(\lambda_i)$ for $i \in S_1$ and Bi$(\tilde{n},\tilde{p}_i)$ for $i \in S_2$, where $S_1 \cup S_2 = \{1,2,\dotsc,n\}$ with $\tilde{p}_i < 1/2$, then
\begin{equation}
d_{TV}(Y,Z) \le \frac{1}{\alpha q} \left(q \sum_{i \in S_1} \lambda_i + \tilde{n} \sum_{i \in S_2} \left(\frac{\tilde{p}_i}{\tilde{q}_i}+q\right) \frac{\tilde{p}_i \tilde{q}_i}{(1-2\tilde{p}_i)^2}\right).\label{pbop}
\end{equation}
\end{itemize}
\end{corollary}

\begin{remarks}
\begin{itemize}
\item[(i)] If $a_{i,l+1} - q a_{i,l} \ge 0$ in Theorem \ref{fmth}, then using the definition of $G_{X_i}$, it is easy to see that 
$$d_{TV}(Y,Z) \le\displaystyle{ \frac{\sigma^2}{\mu}-\frac{q}{p}}.$$
\item[(ii)] The bound given in (\ref{gbop}) is same as the one given in (14), p. 101, of Vellaisamy et. al. \cite{VUC}, which is of constant order. Note that the approach used in proof is more general and easier than approach used in Vellaisamy et. al. \cite{VUC}. Also, it is an improvement over Theorem $2.2$ of Vellaisamy and Upadhye \cite{VN} and comparable to Theorem $1$ of Roos \cite{R3}.
\item[(iii)] If we replace $p_i = p,~i=1,2,\dotsc,n$ in (\ref{gbop}), then bound is exact, as expected.
\item[(iv)] (\ref{pgop}), (\ref{bgop}) and (\ref{pbop}) give bounds for sum of two different types of rvs and can be easily extended for more than two different types of rvs.  
\item[(v)] Instead of multiplying $(1-q z)$ in proof of Theorem \ref{fmth}, we can multiply appropriate function to get the perturbation of some other known distribution and hence the technique used can be generalized. 
\item[(vi)] The bound, in (\ref{pbop}), is not a good bound, as $Y$ has mean greater than variance but in NB variance is bigger than mean, as expected. 
\end{itemize}
\end{remarks}

\subsection{Two-parameter approximation}
Next, we derive the bound between $Y$ and $Z$ by matching first two moments, mean and variance, as
\begin{equation}
\frac{\alpha q}{p} = \mu~~{\rm and}~~ \frac{\alpha q}{p^2} = \sigma^2 \implies p = \frac{\mu}{\sigma^2}~~{\rm and}~~ \alpha = \frac{\mu^2}{\mu_2},\label{mmvariance}
\end{equation}
where $\mu$, $\mu_2$ and $\sigma^2$ are defined as in (\ref{notation}).

\begin{theorem}\label{smth}
Let $X_1,X_2,\dotsc,X_n$ independent rvs with (\ref{assumption}) and $Y=\sum_{i=1}^{n} X_i$, then 
$$d_{TV}\left(Y, Z\right) \le  \frac{1}{\alpha q} \sqrt{\frac{2}{\pi}}\left(\frac{1}{4} + \sum_{i=1}^{n}(1 - d_{TV}(X_i, X_{i+1}))\right)^{-1/2}\sum_{i=1}^{n} \sum_{l=1}^{\infty}l(l-1)\left|a_{i,l+1} - q a_{i,l}\right|,$$
where $\sigma^2 > \mu$ and $Z \sim {\rm NB}(\alpha,p).$
\end{theorem}

\noindent
{\bf Proof}. Using Newton's expansion,  
\begin{equation}
\Delta g(m+j) = \sum_{u=1}^{j-1} \Delta^2 g(m+u) + \Delta g(m+1).\label{ne}
\end{equation}
Substituting (\ref{delta}) and (\ref{ne}) in (\ref{opertilde}), we get
\begin{align}
{\cal A}_Y g(m) &= q\left(\alpha+m\right)g(m+1) - m g(m) + \left[\left(\sum_{i=1}^{n} a_{i,1} - \alpha q \right) + \sum_{i=1}^{n} \sum_{l=1}^{\infty}\left(a_{i,l+1} - a_{i,l}\right)\right]g(m+1)\nonumber\\
&~~+ \sum_{i=1}^{n} \sum_{l=1}^{\infty} l \left(a_{i,l+1} - q a_{i,l}\right) \Delta g(m+1)+ \sum_{i=1}^{n}\sum_{l=1}^{\infty}\sum_{j=1}^{l} \sum_{u=1}^{j-1}\Delta^2 g(m+u)\left(a_{i,l+1} - q a_{i,l}\right),\label{mmoper}
\end{align}
where $q = 1-p$ and $\alpha$, $p$ is defined in (\ref{mmvariance}). Using (\ref{mmvariance}) in (\ref{mmoper}), we obtain
$${\cal A}_Y g(m) = q(\alpha+m)g(m+1) - m g(m) + \sum_{i=1}^{n}\sum_{l=1}^{\infty}\sum_{j=1}^{l} \sum_{u=1}^{j-1}\Delta^2 g(m+u)\left(a_{i,l+1} - q a_{i,l}\right)= {\cal A}_Z g(m) + \tilde{\cal U}_Y g(m).$$
This is a Stein operator of sum of $n$ independent rvs by matching mean and variance with NB rv. Now, taking expectation of $\tilde{\cal U}_Y$ w.r.t $Y$, we have
\begin{align*}
{\mathbb E}[{\tilde{\cal U}_Y g(Y)}] &= \sum_{m=0}^{\infty}\left(\sum_{i=1}^{n}\sum_{l=1}^{\infty}\sum_{j=1}^{l} \sum_{u=1}^{j-1}\Delta^2 g(m+u)\left(a_{i,l+1} - q a_{i,l}\right)\right){\bf P}(Y=m)\\
&=\sum_{m=0}^{\infty}\sum_{i=1}^{n}\sum_{l=1}^{\infty}\sum_{j=1}^{l} \sum_{u=1}^{j-1} \Delta g(m+u) \left[{\bf P}(Y=m-1) - {\bf P}(Y=m)\right] \left(a_{i,l+1} - q a_{i,l}\right).
\end{align*}
Therefore, for $g \in {\cal G}_Z \cap {\cal G}_Y$, we have
$$|{\mathbb E}[{\tilde{\cal U}_Y g(Y)}]| \le  d_{TV}(Y, Y+1)\|\Delta g\|\sum_{i=1}^{n} \sum_{l=1}^{\infty} l(l-1)\left|a_{i,l+1} - q a_{i,l}\right|.$$
Then, the proof is follows by using (\ref{mros}) and (\ref{bound}).\qed

\begin{corollary}\label{cor2}
Given $Y=\sum_{i=1}^{n} X_i$, let us choose $X_i$ different type of distribution, then we have the following bounds
\begin{itemize}
\item[{(i)}]Let $X_i$ follow Ge$(p_i),~i=1,2,\dotsc,n$ with $q_i = (1-p_i) < 1/2$, then
\begin{equation}
d_{TV}(Y,Z) \le \frac{2}{\mu} \sqrt{\frac{2}{\pi}} \left(\sum_{i=1}^{n} q_i + \frac{1}{4}\right)^{-1/2} \sum_{i=1}^{n} \left|\frac{1}{p}-\frac{1}{p_i}\right| \left(\frac{q_i}{p_i}\right)^2.\label{gbtp}
\end{equation}
\item[{(ii)}] Let $X_i$ follow Po$(\lambda_i)$ for $i \in S_1$ and Ge$(p_i)$ for $i \in S_2$, where $S_1 \cup S_2 = \{1,2,\dotsc,n\}$ and $q_i < 1/2$, then
\begin{equation}
d_{TV}(Y,Z) \le \frac{2}{\alpha q} \sqrt{\frac{2}{\pi}} \left(\frac{1}{4} + \sum_{i\in S_1} \left(1- \frac{e^{-\lambda_i} \lambda_i^{\floor{\lambda_i}}}{\floor{\lambda_i}!}\right)+\sum_{i \in S_2}q_i\right)^{-1/2} \left(\sum_{i \in S_2} |p-p_i| \frac{q_i^2}{p_i^3}\right).\label{pgtp}
\end{equation}
\item[{(iii)}] Let $X_i$ follow Bi$(\tilde{n},\tilde{p}_i)$ for $i \in S_1$ and Ge$(p_i)$ for $i \in S_2$, where $S_1 \cup S_2 = \{1,2,\dotsc,n\}$ with $q_i, \tilde{q}_i < 1/2$, then
\begin{equation}
d_{TV}(Y,Z) \le \frac{1}{\alpha q}\sqrt{\frac{2}{\pi}} \frac{\displaystyle{\left(\tilde{n} \sum_{i \in S_1}\left(\frac{\tilde{p}_i}{\tilde{q}_i}+q\right) \frac{\tilde{p}_i^2 \tilde{q}_i}{(1-2\tilde{q}_i)^3}+2 \sum_{i \in S_2} |p-p_i|\frac{q_i^2}{p_i^3}\right)}}{\displaystyle{\small{\left(\frac{1}{4} + \sum_{i\in S_1} \left(1- {\tilde{n} \choose \floor{(\tilde{n}+1)\tilde{p}_i}}\tilde{p}_i^{\floor{(\tilde{n}+1)\tilde{p}_i}} (1-\tilde{p}_i)^{\tilde{n}-\floor{(\tilde{n}+1)\tilde{p}_i}}+\frac{\tilde{p}_i^{\tilde n}}{2}\right)+\sum_{i \in S_2}q_i\right)^{1/2} }}},\label{bgtp}
\end{equation}
where $\displaystyle{\tilde{n} \sum_{i \in S_1} \tilde{p}_i^2<\sum_{i \in S_2} \frac{q_i^2}{p_i^2}}$.
\end{itemize} 
\end{corollary}

\begin{remarks}
\begin{itemize}
\item[(i)] If $p_i = p,i=1,2,\dotsc,n$ in (\ref{gbtp}), then bound is exact, as expected.
\item[(ii)] The bound in (\ref{gbtp}) is improvement by constant over Corollary $4.1$ of Vellaisamy et. al. \cite{VUC}, which is of order $O(n^{-1/2})$.
\item[(iii)] (\ref{pgtp}) and (\ref{bgtp}) give bound for two different types of rvs, where variance is greater than mean. Also, this can be extended for more than two random variables.
\end{itemize}
\end{remarks}

\subsection{Three-parameter approximation} 
As mentioned in Section $2$, NB distribution can be described using two parameters. Therefore, for three parameter approximation, we use convolution of one parameter distribution, namely geometric, with NB. Convolution of Poisson with NB is studied by Vellaisamy et. al. \cite{VUC} and has improved the accuracy of approximation with respect to NB or Poisson approximation. Therefore, we choose geometric distribution as it has behavior similar to NB distribution. \\
Next, we derive a Stein operator for convolution of NB and Geometric. Recall that $Z \sim {\rm NB}(\alpha,p)$ and $W \sim Ge({\hat{p}})$, then the PGF of $Z$ and $W$ is given by $M_Z(z) = p^\alpha/(1-q z)^{\alpha}$ and $M_W(z) = \hat{p}/(1 - \hat{q} z)$ respectively. Also, $V = W + Z$, then the PGF of $V$ is  $M_V(z) = p^{\alpha} {\hat{p}} /\left((1-q z)^{\alpha}(1 - {\hat{q}} z)\right)$. Differentiating with respect to $z$ and multiplying by ($1-qz$), we get 
$$\sum_{m=0}^{\infty} (m+1) p^\prime_{m+1} z^m - q \sum_{m=0}^{\infty} m p^\prime_m z^m  =  \alpha q \sum_{m=0}^{\infty}p^\prime_m z^m + \sum_{m=0}^{\infty}\left(\sum_{l=0}^{m} p^\prime_l {{\hat{q}}}^{m-l+1} \right)z^m - q \sum_{m=0}^{\infty} \left(\sum_{l=0}^{m-1}p^\prime_l {{\hat{q}}}^{m-l}\right)z^m,$$
where $p^\prime_m = {\bf P}(V =m)$ be the PMF of $V$. Comparing the coefficient of $z^m$, we have
$$\alpha q p^\prime_m + q m p^\prime_m - (m+1)p^\prime_{m+1} = q \sum_{l=0}^{m-1}p^\prime_l {{\hat{q}}}^{m-l} - \sum_{l=0}^{m} p^\prime_l {{\hat{q}}}^{m-l+1}= \sum_{l=0}^{m}(q - {{\hat{q}}}) p^\prime_l {{\hat{q}}}^{m-l} - q p^\prime_m,$$
This can be written as
$$q(\alpha +1 + m) p^\prime_m - (m+1)p^\prime_{m+1} + \sum_{l=0}^{m}({{\hat{q}}} - q) p^\prime_l {{\hat{q}}}^{m-l} = 0.$$
For $g \in {\cal G}_V$, defined in (\ref{gx}), we have
$$\sum_{m=0}^{\infty} g(m+1) \left[ (\alpha+1) q p^\prime_m + q m p^\prime_m - (m+1)p^\prime_{m+1}+\sum_{l=0}^{m}({{\hat{q}}} - q) p^\prime_l{{\hat{q}}}^{m-l}\right] = 0$$
Hence,
$${\mathbb E}[{\cal A}_V g(V)] = \sum_{m=0}^{\infty} \left[q (\alpha+1 + m) g(m+1) - m g(m) + ({\hat{q}} -q)\sum_{l=0}^{\infty} g(m+l+1) {{\hat{q}}}^l\right]p^\prime_m = 0,$$
where ${\cal A}_V g(m) = q (\alpha +1 + m) g(m+1) - m g(m) + ({\hat{q}}-q)\sum_{l=0}^{\infty} g(m+l+1) {{\hat{q}}}^l$ is a Stein operator for $V$, which is a perturbation of NB$(\alpha+1, p)$. Using (\ref{delta}), Stein operator can be written as
$${\cal A}_V g(m) = q \left[\left(\alpha + 1 + \frac{\hat{q} - q}{q \hat{p}}\right) + m\right] g(m+1) - m g(m) + \left(\frac{\hat{q} - q}{\hat{p}}\right) \sum_{j=1}^{\infty} \Delta g(m+j) {\hat{q}}^l= {\hat{\cal A}_Z}g(m) + \hat{\cal U}_V g(m),$$
where $\hat{\cal A}_Z$ is a Stein operator for NB($r,p$) with $r=\left(\displaystyle{\alpha + 1 + \frac{\hat{q} - q}{q \hat{p}}}\right)$. Then
\begin{equation}
{\mathbb E}[\hat{\cal U}_V g(V)] \le \|\Delta g\|\times \left|\hat{q}-q\right| \frac{\hat{q} }{\hat{p}^2}.\label{pert}
\end{equation}
Next, we match the first three moments of $V$ and $Y$, we have 
\begin{align}
\frac{\alpha q}{p} + \frac{{\hat{q}}}{{\hat{p}}} &= \mu,  \quad \frac{\alpha q^2}{p^2} + \frac{{{\hat{q}}}^2}{{{\hat{p}}}^2} = {\mu_2}\quad {\rm and} \quad \frac{\alpha q^3}{p^3} + \frac{{{\hat{q}}}^3}{{{\hat{p}}}^3} = \frac{\mu_3}{2}, \label{tpmatching}
\end{align}
where $\mu$, $\mu_2$ and $\mu_3$ defined as in (\ref{notation}). Therefore, the choice of parameters is
\begin{equation}
{\hat{p}} =  \frac{3 \mu}{\left(3 \mu + \eta\right)}, \quad \alpha = \left(\mu- \frac{\eta}{3 \mu}\right)^2\Bigr/ \left(\mu_2 - \frac{\eta^2}{9 \mu^2}\right)~{\rm and}~p  = \left(\mu -\frac{\eta}{3 \mu}\right)\Bigr/ \left(\sigma^2 - \frac{\eta}{3 \mu}\left(\frac{\eta}{3 \mu}+1\right)\right),\label{tpmatching2}
\end{equation}
where $\eta $ as defined in (\ref{eta}). Now, we are obtain bound for $V$ approximation to $Y$ by matching first three moments.

\begin{theorem}\label{thmth}
Let $X_1,X_2,\dotsc,X_n$ are independent rvs with (\ref{assumption}) and $\sigma^2 > \mu$, then
\begin{align*}
d_{TV}\left(Y, V \right) &\le  \frac{16}{\Psi \times \left(rq-|\hat{q}-q|\displaystyle{\frac{\hat{q}}{\hat{p}^2}}\right)}  \left(\sum_{i=1}^{k}\sum_{l=1}^{\infty} \frac{l(l-1)(l-2)}{6}\left|a_{i_{l+1}}- q a_{i_l}\right| + \frac{|{\hat{q}} - q|{\hat{q}}^3}{{\hat{p}}^4}  \right),
\end{align*}
where $\displaystyle{\Psi = \sum_{i=1}^{n} \xi_i,~ \xi_i =\min_{1 \le i \le n} \left(\frac{1}{2},~1-d_{TV}(X_i,X_i+1)\right)} $ and $rq>|\hat{q}-q|\displaystyle{\frac{\hat{q}}{\hat{p}^2}}$.
\end{theorem}

\noindent
{\bf Proof}. Now, for $g \in {\cal G}_Y$, we introduce a parameter $\hat{q}$ and modify Stein operator of $Y$ in (\ref{mmoper}) as follows
\begin{align*}
{\cal A}_Y g(m) &= q(\alpha+1+m)g(m+1) - m g(m) + ({\hat{q}} - q)\sum_{l=0}^{\infty}g(m+l+1){{\hat{q}}}^l\\
&~+ \left[\hspace{-0.1cm}\left(\sum_{i=1}^{n} a_{i,1} \hspace{-0.1cm}-\hspace{-0.1cm} \alpha q \hspace{-0.1cm}- q\right)\hspace{-0.1cm} +\hspace{-0.1cm} \sum_{i=1}^{n} \sum_{l=1}^{\infty}\hspace{-0.1cm}\left(a_{i,l+1} \hspace{-0.1cm}-\hspace{-0.1cm} q a_{i,l}\right)\right]\hspace{-0.1cm}g(m+1)\hspace{-0.1cm} +\hspace{-0.1cm} \sum_{i=1}^{n} \sum_{l=1}^{\infty} l\hspace{-0.1cm} \left(a_{i,l+1} \hspace{-0.1cm}-\hspace{-0.1cm} q a_{i,l}\right) \hspace{-0.1cm}\Delta g(m+1)\\
&~+ \sum_{i=1}^{n}\sum_{l=1}^{\infty}\sum_{j=1}^{l} \sum_{u=1}^{j-1}\Delta^2 g(m+u)\left(a_{i,l+1} - q a_{i,l}\right) - ({\hat{q}} - q)\sum_{l=0}^{\infty}g(m+l+1){{\hat{q}}}^l\\
&= \hat{\cal A}_V g(m) +\hat{\cal U}_Y g(m),
\end{align*}
where $q = 1-p$ and $\hat{q} = 1 - \hat{p}$. Also, $\alpha$, $p$ and $\hat{p}$ is defined in (\ref{tpmatching2}). Again, from Newton's expansion, we have
\begin{equation}
\Delta^2 g(m+u) = \sum_{v=1}^{u-1} \Delta^3 g(m+v) + \Delta^2 g(m+1).\label{ne2}
\end{equation}
Substituting (\ref{delta}), (\ref{ne}) and (\ref{ne2}) in $\hat{\cal U}_Y$ then using (\ref{tpmatching}), we get
\begin{equation}
\hat{\cal U}_Y g(m) = \sum_{i=1}^{n}\sum_{l=1}^{\infty} \sum_{j=1}^{l} \sum_{u=1}^{j-1} \sum_{v=1}^{u-1} \Delta^3 g (m+v) \left(a_{i,l+1} - q a_{i,l}\right)- \frac{({\hat{q}} - q){{\hat{q}}}^2}{{\hat{p}}^3} \sum_{v=1}^{\infty} \Delta^3 g(m+v){{\hat{q}}}^v.\label{ppop}
\end{equation}
Taking expectation w.r.t $Y$, we have
\begin{align*}
{\mathbb E} (\hat{\cal U}_Y g(Y)) &= \sum_{m=0}^{\infty} \left(\sum_{i=1}^{n}\sum_{l=1}^{\infty} \sum_{j=1}^{l} \sum_{u=1}^{j-1} \sum_{v=1}^{u-1} \Delta^3 g(m+v) \left(a_{i,l+1}- q a_{i,l}\right)- \frac{({\hat{q}} - q){{\hat{q}}}^2}{{\hat{p}}^3} \sum_{v=1}^{\infty} \Delta^3 g(m+v){{\hat{q}}}^v \right) p_m\\
&= \sum_{m=0}^{\infty}  \left(\sum_{i=1}^{n}\sum_{l=1}^{\infty} \sum_{j=1}^{l} \sum_{u=1}^{j-1} \sum_{v=1}^{u-1} \Delta g(m+v) [p_{m-2} - 2 p_{m-1} + p_m]  \left(a_{i,l+1}- q a_{i,l}\right)\right.\\
&~~~~~~~~~~~~~~~~~~~~~~~~~~~~~~~~~~~~~~~~~~~~~~~~~~~\left.- \frac{({\hat{q}} - q){{\hat{q}}}^2}{{\hat{p}}^3} \sum_{v=1}^{\infty}\Delta g(m+v) [p_{m-2} - 2 p_{m-1} + p_m] {{\hat{q}}}^v \right)
\end{align*}
Hence, 
\begin{align}
|{\mathbb E} (\hat{\cal U}_Y g(Y))| &\le \|\Delta g\|\left(\sum_{i=1}^{n}\sum_{l=1}^{\infty}\frac{l(l\hspace{-0.1cm}-\hspace{-0.1cm}1)(l\hspace{-0.1cm}-\hspace{-0.1cm}2)}{6} \left|a_{i,l+1}\hspace{-0.1cm}-\hspace{-0.1cm} q a_{i,l}\right| + \frac{|{\hat{q}} \hspace{-0.1cm}-\hspace{-0.1cm} q|{\hat{q}}^3}{{\hat{p}}^4}  \right)\sum_{m=0}^{\infty} |p_{m-2} \hspace{-0.1cm}-\hspace{-0.1cm} 2 p_{m-1}\hspace{-0.1cm} + \hspace{-0.1cm}p_m|. \label{sop}
\end{align}
From ($4.9$) of Barbour and \v{C}ekanavi\v{c}ius \cite{BC}, we have
$$\sum_{m=0}^{\infty} |p_{m-2} - 2 p_{m-1} + p_m| \le \frac{16}{\Psi},$$
where $\Psi = \sum_{i=1}^{k} \xi_i$ and $\xi_i = \min_{1 \le i \le n} \left\{\frac{1}{2}, d_{TV}\left(X_i,X_i+1\right)\right\}$.\\ 
Hence, from (\ref{perturbation}) and (\ref{pert}), the proof follows.\qed

\begin{corollary}\label{cor3}
Let $Y=\sum_{i=1}^{n} X_i$ such that $X_i$ be the different types of distribution, then for $rq>|\hat{q}-q|\displaystyle{\frac{\hat{q}}{\hat{p}^2}}$, we have the following bounds
\begin{itemize}
\item[{(i)}] Let $X_i$ follow Ge$(p_i),~i=1,2,\dotsc,n$ with $q_i = (1-p_i) < 1/2$, then
\begin{equation}
d_{TV}(Y,V) \le \frac{16 \left(\sum_{i=1}^{n}q_i\right)^{-1}}{p\times\left(rq-|\hat{q}-q|\displaystyle{\frac{\hat{q}}{\hat{p}^2}}\right)} \left(\sum_{i=1}^{n} \left|\frac{1}{p}-\frac{1}{p_i}\right| \left(\frac{q_i}{p_i}\right)^3 + \left|\frac{1}{p} - \frac{1}{\hat{p}}\right| \left(\frac{\hat{q}}{\hat{p}}\right)^3\right).\label{gbthp}
\end{equation} 
\item[{(ii)}] Let $X_i$ follow Po$(\lambda_i)$ for $i \in S_1$ and Ge$(p_i)$ for $i \in S_2$, where $S_1 \cup S_2 = \{1,2,\dotsc,n\}$ and $q_i < 1/2$, then
\begin{align}
d_{TV}\left(Y, V \right) &\le  \frac{16 \Psi^{-1}}{p \times \left(rq-|\hat{q}-q|\displaystyle{\frac{\hat{q}}{\hat{p}^2}}\right)}  \left(\sum_{i \in S_2}\left|\frac{1}{p}-\frac{1}{p_i}\right|\left(\frac{q_i}{p_i}\right)^3 + \left|\frac{1}{\hat{p}} -\frac{1}{p}\right| \left(\frac{{\hat{q}}}{{\hat{p}}}\right)^3  \right),\label{pgthp}
\end{align}
\item[{(iii)}] Let $X_i$ follow Bi$(\tilde{n},\tilde{p}_i)$ for $i \in S_1$ and Ge$(p_i)$ for $i \in S_2$, where $S_1 \cup S_2 = \{1,2,\dotsc,n\}$ with $q_i, \tilde{p}_i < 1/2$, then
\begin{align}
d_{TV}\left(Y, V \right) &\le  \frac{\displaystyle{16\left(\frac{\tilde n}{p} \sum_{i \in S_1} \left(\frac{\tilde{p}_i}{\tilde{q}_i} + q\right)\frac{\tilde{p}_i^3 \tilde{q}_i}{(1-2\tilde{p}_i)^4} + \sum_{i \in S_2}\left|\frac{1}{p}-\frac{1}{p_i}\right|\left(\frac{q_i}{p_i}\right)^3 + \left|\frac{1}{\hat{p}} - \frac{1}{p}\right|\left(\frac{\hat{q}}{\hat{p}}\right)^3  \right)}}{\Psi \times p\left(rq-|\hat{q}-q|\displaystyle{\frac{\hat{q}}{\hat{p}^2}}\right)},\label{bgthp}
\end{align}
where $\displaystyle{\tilde{n} \sum_{i \in S_1} \tilde{p}_i^2<\sum_{i \in S_2} \frac{q_i^2}{p_i^2}}$.
                                                                                                                                                                                                                                                                                                                                                                                                                                                                              \end{itemize}
\end{corollary}

\begin{remarks}
\begin{itemize}
\item[(i)] If $p_i = p,i=1,2,\dotsc,n$ and $\hat{p} = p$ in (\ref{gbthp}), then bound is exact, as expected.
\item[(ii)] The bound in Theorem \ref{thmth} is of order $O(n^{-1})$, which is improvement over one and two parameter approximation.
\item[(iii)] (\ref{gbthp}) and (\ref{pgthp}) give bounds for all are one type rvs and two different type of rvs. Also, this can extend to more than two types of rvs.
\item[(iv)] We can not obtain bound for sums of binomial and Poisson rvs for two and three parameter approximation because mean is greater than variance. So, the choice of parameters is inadmissible. 
\end{itemize}
\end{remarks}

\section{An Application}\label{app}
In this section, we demonstrate an application of our approximation results to obtain bound between NB and a function of waiting time for binomial distribution of order $(k_1,k_2)$ (see Huang and Tsai \cite{HT}).\\ 
Let {\bf S} denote success and {\bf F} failure, with success probability $\bar{p}$, in a sequence of independent Bernoulli trials. If $k_1$ consecutive {\bf F}s followed by $k_2$ consecutive {\bf S}s, i.e.,
$$\dotsc \overbrace{{\bf F} \dotsc{\bf F}}^{k_1}\overbrace{{\bf S}\dotsc{\bf S}}^{k_2}\dotsc,$$ 
occurred then it is called $(k_1,k_2)$-event, where $(k_1,k_2)$ is a pair of nonnegative integers, including $0$, excluding $(0,0)$. Also, let $\tilde{N}(n;~k_1,~k_2)$ be number of occurrences of $(k_1,k_2)$-events in $n$ trials. The distribution of $\tilde{N}(n;~k_1,~k_2)$, denoted by $p_{\cdot, n}$, is called the binomial distribution of order $(k_1,k_2)$. $p_{\cdot,n}$ defined in Lemma $1$ of Huang and Tsai \cite{HT} as follows:
\begin{lemma}\leavevmode
\begin{itemize}
\item[(i)]  
$p_{x,n} = \left\{
\begin{array}{l l}
             0  \quad & \text{if $n < k_1+k_2,~x > 0$};\\
             1 \quad &  \text{if $n < k_1+k_2,~x=0$};\\
             \bar{q}^{k_1} \bar{p}^{k_2} \quad & \text{if $n=k_1+k_2,~x=1$};\\
             1-\bar{q}^{k_1} \bar{p}^{k_2} \quad & \text{if $n=k_1+k_2,~x=0$}.
\end{array}\right.$
\item[(ii)] $\displaystyle{p_{0,n} = p_{0,n-1}-\bar{q}^{k_1} \bar{p}^{k_2}~p_{0,n-k_1-k_2}}$.
\item[(iii)] $\displaystyle{p_{x,n}= \sum_{j=0}^{n-k_1-k_2}\bar{q}^{k_1} \bar{p}^{k_2}p_{x-1,n}~p_{0,n-k_1-k_2}}$.
\item[(iv)] $\displaystyle{p_{x,n+1}= p_{x,n}+\bar{q}^{k_1} \bar{p}^{k_2}\left[p_{x-1,n-k_1-k_2+1}-p_{x,n-k_1-k_2+1}\right] \quad {\rm for}~ n \ge k_1+k_2,~1 \le x \le \floor{\frac{n}{k_1+k_2}}}$,\\
where $\floor{a}$ denote the greatest integer not exceeding $a$.
\end{itemize}
\end{lemma}
\noindent
Next, let $\tilde{T}_n$ denote the waiting time for $n$th occurrence of $(k_1,~k_2)$-event. Then
$$\tilde{T}_n = T_1+T_2+\dotsb+T_n,$$
where $T_j$ is $k:=k_1+k_2$ plus the number of trials between the $(j-1)^{th}$ and $j^{th}$ occurrence of $(k_1,~k_2)$-event. $T_j$'s are independent and identically distributed (i.i.d.) with i.i.d. copy $T$ having PMF
\[
{\bf P}(T = n) = \left\{
\begin{array}{l l}
           0 \quad & \text{$n < k$};\\
           a(\bar{p}) \quad & \text{$n=k$};\\
           a(\bar{p}) p_{0,n-k}  \quad &\text{$n>k$},
\end{array}\right. 
\]  
where $a(\bar{p}) = {\bar{q}}^{k_1} {\bar{p}}^{k_2}$. \\
Define $\hat{T}_j = T_j-k$, for $j=1,2,\dotsc,n$. Therefore, $\hat{T}_j$ is the number of trials between $(j-1)^{th}$ and $j^{th}$ occurrence of $(k_1,k_2)$ event. Suppose $\hat{T}$ be the i.i.d. copy of $\hat{T}_j$. Then
\[
{\bf P}(\hat{T}=n) = {\bf P}(T = n+k) = \left\{
\begin{array}{l l}
           0  \quad & \text{$n < 0$};\\
           a(\bar{p}) \quad & \text{$n=0$};\\
           a(\bar{p}) p_{0,n}  \quad &\text{$n>0$}.
\end{array}\right. 
\]  
Let $M_{\hat{T}}(t)$ be the PGF of $\hat{T}$. Then, it can be easily seen that
$$M_{\hat{T}}(z) =\frac{a(\bar{p})}{1-z+a(\bar{p})z^k}$$
(see Huang and Tsai \cite{HT}, p.p. 128-129, for details). Define $\acute{T}$ as, the number of failures before $n^{th}$ occurrence of $(k_1,k_2)$-event,
\begin{equation}
\acute{T} = \hat{T}_1 + \hat{T}_2 + \dotsb + \hat{T}_n.\label{fevent}
\end{equation}
Then, the PGF of $\acute{T}$ is
\begin{equation}
M_{\acute{T}}(z) = \left(\frac{a(\bar{p})}{1-z+a(\bar{p})z^k}\right)^n. \label{wtpgf}
\end{equation}
Also, define
\begin{equation}
b_{m,\bar{p}} :=\sum_{l=0}^{\floor{m/k}}(-1)^l {m-l(k-1) \choose l} a(\bar{p})^l.\label{nots}
\end{equation}
For more details of $(k_1,k_2)$ distribution, we refer the readers to Balakrishnan and Koutras \cite{BALA}, Dafnis et. al. \cite{DAP} and Huang and Tsai \cite{HT}.

\subsection{One-parameter approximation}
First, we derive bound between $\acute{T}$ and $Z$ by matching first moment as follows: 
$$\frac{\alpha q}{p} = n \left(\frac{1-k a(\bar{p})}{a(\bar{p})}\right).$$
Here, matching can be done in two ways:
\begin{enumerate}
\item Let $p$ be fixed, of our choice, and $\alpha = n p (1-k a(\bar{p}))/q a(\bar{p}) $.
\item Let $\alpha$ be fixed and $p=\alpha a(\bar{p})/(\alpha a(\bar{p})+n(1-ka(\bar{p})))$.
\end{enumerate}
For one-parameter approximation, we fixed $\alpha q = n$ and 
\begin{equation}
p=\frac{a(\bar{p})}{(1-ka(\bar{p}))}.\label{Mmean2}
\end{equation}
\begin{theorem}\label{aoop}
Let $\acute{T}$ be defined in (\ref{fevent}), then
\begin{equation}
d_{TV}(\acute{T},Z) \le\left(1-k a(\bar{p})\right)\hspace{-0.1cm}\sum_{l=1}^{\infty} l \left|b_{l,\bar{p}}- q b_{l-1,\bar{p}}\right|\hspace{-0.1cm}+\hspace{-0.1cm} k(k-1) a(\bar{p})+k a(\bar{p}) \left(\hspace{-0.1cm}\sum_{l=1}^{\infty} \left|b_{l,\bar{p}} - q b_{l-1,\bar{p}}\right| +1\hspace{-0.1cm} \right),\label{aop}
\end{equation}
where $Z \sim {\rm NB}(\alpha, p)$, and $b_{l,\bar{p}}$ as defined in (\ref{nots}).
\end{theorem}
\begin{remarks}
\begin{itemize}
\item[(i)] The bound in Theorem \ref{aoop} is of constant order and can be calculated for different values of $(k_1,k_2)$ and $\bar{p}$.
\item[(ii)] In Table $1$, the bound for various values of $(k_1,k_2)$ and $\bar{p}$ is calculated by taking $l$ up to  first $3000$ terms and neglecting the remainder, as the values are too small. Also, we can observe the pattern that, as the value of $\bar{p}$ decreases the bound decreases, which is consistent with NB convergence to Poisson.
\end{itemize}
\end{remarks}

\begin{table}[htbp]
\resizebox{0.83\textwidth}{!}{\begin{minipage}{\textwidth}
\caption{One-Parameter Approximation}
\begin{center}
\setlength\extrarowheight{1pt}
\begin{tabular}{c c}
\begin{tabular}{|c|c|c|c|}
\hline
$(k_1,k_2)$ & $\bar{p}=1/4$ & $\bar{p}=1/8$ & $\bar{p}=1/16$\\
\hline
$(1,4)$ & $1.05816$ & $0.141903$ & $0.0013244$\\
\hline
$(1,5)$ & $0.672985$ & $0.00424482$ & $0.0000438252$\\
\hline
$(1,6)$ & $0.116891$ & $0.000260107$ & $3.53445 \times 10^{-6}$\\
\hline
$(1,7)$ & $0.0124531$ & $0.0000358306$ & $2.93422 \times 10^{-7}$\\
\hline
$(1,8)$ & $0.00181219$  & $5.64488 \times 10^{-6}$ & $2.35743 \times 10^{-8}$\\
\hline
$(1,9)$ & $0.000422883$ & $8.80291 \times 10^{-7}$ & $1.84173 \times 10^{-9}$\\
\hline
$(2,4)$ & $1.05743$ & $0.117485$ & $0.00139172$\\
\hline
$(2,5)$ & $0.521417$ & $0.00381735$ & $0.0000559635$\\
\hline
$(2,6)$ & $0.0762276$ & $0.000283511$ & $4.41285 \times 10^{-6}$\\
\hline
$(2,7)$ & $0.00866093$ & $0.0000400281$ & $3.5366 \times 10^{-7}$\\
\hline
\end{tabular}
&
\begin{tabular}{|c|c|c|c|}
\hline
$(2,8)$ & $~0.00148472~$ & $6.17008 \times 10^{-6}$ & $2.76261 \times 10^{-8}$\\
\hline
$(3,4)$ & $1.03096$ & $0.0975824$ & $0.00148602$\\
\hline
$(3,5)$ & $0.382432$ & $0.00352341$ & $0.0000687831$\\
\hline
$(3,6)$ &  $0.0500174$ & $0.000305253$ & $5.31503 \times 10^{-6}$\\
\hline
$(3,7)$ & $~0.00631033~$ & $0.0000435841$ & $4.14431 \times 10^{-7}$\\
\hline
$(4,4)$& $0.961538$ & $0.0814895$ & $0.00160097$\\
\hline
$(4,5)$ & $0.26916$ & $0.00332182$ & $0.0000820025$\\
\hline
$(4,6)$ & $~0.0334377~$ & $~0.000324287~$ & $~6.22537 \times 10^{-6}~$\\
\hline
$(5,4)$ & $0.844592$ & $0.0685582$ & $0.00173112$\\
\hline
$~~(5,5)~~$ & $~0.184463~$ & $~0.00318198~$ & $~0.0000953822~$\\
\hline
$~(6,4)~$ & $~~0.69535~~$ & $~~0.0582122~~$ & $~~~0.0018718~~~$\\
\hline
\end{tabular}
\end{tabular}
\end{center}
\end{minipage}}
\label{tab:multicol}
\end{table}

\subsection{Two-parameter approximation}
Next, we derive bound between $\acute{T}$ and $Z$ by matching mean and variance as
\begin{equation}
\frac{\alpha q}{ p} = \frac{n(1-ka(\bar{p}))}{a(\bar{p})} \quad{\rm and}\quad \frac{\alpha q}{p^2} = n \left(\frac{1-(2k-1)a(\bar{p})}{a(\bar{p})^2}\right).\label{mvar}
\end{equation}
This leads to the following choice of parameters
\begin{equation}
p = \frac{(1-ka(\bar{p}))a(\bar{p})}{1-(2k-1)a(\bar{p}) } \quad {\rm and } \quad \alpha = \frac{n (1-k a(\bar{p}))^2}{1-2k a(\bar{p})+ka(\bar{p})^2}.\label{alphap}
\end{equation}
\begin{theorem}\label{atop}
Let $\acute{T}$ be defined in (\ref{fevent}) with $1-2k a(\bar{p})+ka(\bar{p})^2 > 0$, then
\begin{equation}
\begin{aligned}
d_{TV}(\acute{T},Z) \le \frac{n}{\alpha q} \sqrt{\frac{2}{\pi}} \left(\frac{1}{4}+ n\left(1- \frac{a(\bar{p})}{2}(1+a(\bar{p}))\right)\right)^{-1/2} \left[\sum_{l=1}^{\infty} \frac{l(l-1)}{2}\left|b_{l,\bar{p}}-q b_{l-1,\bar{p}}\right|\right.\\
\left. +\frac{k(k-1)(k-2)}{2}a(\bar{p})+k a(\bar{p})\sum_{l=k}^{\infty} \frac{l(l-1)}{2}\left|b_{l-k+1,\bar{p}}- q b_{l-k,\bar{p}}\right|\right], \label{atp}
\end{aligned}
\end{equation}
where $Z \sim {\rm NB}(\alpha,p)$ and $b_{l,\bar{p}}$ as defined in (\ref{nots}).
\end{theorem}
\begin{remarks}
\begin{itemize}
\item[(i)] The bound in Theorem \ref{atop} is of order $O(n^{-1/2})$. Therefore, as $n$ increases the bound decreases.
\item[(ii)] It is easy to see that the bound in two parameter approximation is better then one parameter approximation (see Table $1$ and Table $2$), as expected.
\item[(iii)] In Table $2$, the bound for various values of $(k_1,k_2)$ and $\bar{p}$ is calculated by taking $l$ up to  first $3000$ terms and neglecting the remainder, as the values are too small. Also, we can observe the pattern that, as the value of $\bar{p}$ decreases the bound decreases, which is consistent with NB convergence to Poisson.
\end{itemize}
\end{remarks}

\begin{table}[htbp]
\resizebox{0.9\textwidth}{!}{\begin{minipage}{\textwidth}
\caption{Two-Parameter Approximation}
\setlength\extrarowheight{1pt}
\begin{center}
\begin{tabular}{|c|c|c|c|c|c|c|}
\hline  
\multicolumn{1}{|c|}{}  & \multicolumn{2}{|c|}{$\bar{p}=1/4$} & \multicolumn{2}{|c|}{$\bar{p}=1/8$} & \multicolumn{2}{|c|}{$\bar{p}=1/16$}\\
\hline
$(k_1,k_2)$ & $n=50$ & $n=100$ & $n=50$ & $n=100$ & $n=50$ & $n=100$ \\
\hline
$(1,4)$&$1.1293$&$0.799532$&$0.0318133$&$0.0225235$&$0.0000787781$&$0.0000557739$\\
\hline
$(1,5)$&$0.645954$&$0.457328$&$0.00037684$&$0.000266798$&$8.05622 \times 10^{-6}$&$5.70371 \times 10^{-6}$\\
\hline
$(1,6)$&$0.046521$&$0.0329363$&$0.0000529935$&$0.0000375187$&$8.80547 \times 10^{-7}$&$6.23417 \times 10^{-7}$\\
\hline
$(1,7)$&$0.00238398$&$0.00168783$&$0.0000105207$&$7.44852 \times 10^{-6}$&$8.80545 \times 10^{-8}$&$6.23415 \times 10^{-8}$\\
\hline
$(1,8)$&$0.000459553$&$0.000325358$&$1.97243 \times 10^{-6}$&$1.39645 \times 10^{-6}$&$8.25511 \times 10^{-9}$&$5.84452\times 10^{-9}$\\
\hline
$(1,9)$&$0.000155114$&$0.000109819$&$3.52218 \times 10^{-7}$&$2.49366 \times 10^{-7}$&$7.37063 \times 10^{-10}$&$5.21832 \times 10^{-10}$\\
\hline
$(2,4)$&$1.64106$&$1.16185$&$0.0345507$&$0.0244615$&$0.00013857$&$0.0000981058$\\
\hline
$(2,5)$&$0.553509$&$0.391878$&$0.000497124$&$0.000351958$&$0.0000132146$&$9.35579 \times 10^{-6}$\\
\hline
$(2,6)$&$0.0305958$&$0.0216614$&$0.0000739922$&$0.0000523856$&$1.32082 \times 10^{-6}$&$9.35125 \times 10^{-7}$\\
\hline
$(2,7)$&$0.00198307$&$0.00140399$&$0.0000138079$&$9.77585 \times 10^{-6}$&$1.23827 \times 10^{-7}$&$8.76677 \times 10^{-8}$\\
\hline
$(2,8)$&$0.000477953$&$0.000338385$&$2.46553 \times 10^{-6}$&$1.74557 \times 10^{-6}$&$1.10559 \times 10^{-8}$&$7.82748 \times 10^{-9}$\\
\hline
$(3,4)$&$2.09053$&$1.48007$&$0.0350269$&$0.0247986$&$0.00021872$&$0.000154851$\\
\hline
$(3,5)$&$0.417545$&$0.295617$&$0.000631458$&$0.000447065$&$0.0000198194$&$0.0000140319$\\
\hline
$(3,6)$&$0.0197698$&$0.0139968$&$0.0000969564$&$0.000068644$&$1.8574 \times 10^{-6}$&$1.31502 \times 10^{-6}$\\
\hline
$(3,7)$&$0.0017595$&$0.00124571$&$0.0000172595$&$0.0000122196$&$6.5839 \times 10^{-7}$&$1.17412 \times 10^{-7}$\\
\hline
$(4,4)$&$2.30166$&$1.62955$&$0.0339291$&$0.0240214$&$0.000319874$&$0.000226467$\\
\hline
$(4,5)$&$0.286267$&$0.202673$&$0.00077711$&$0.000550184$&$0.0000278687$&$0.0000197307$\\
\hline
$(4,6)$&$0.013037$&$0.00923002$&$0.000121071$&$0.0000857165$&$2.48759 \times 10^{-6}$&$1.76119 \times 10^{-6}$\\
\hline
$(5,4)$&$2.18446$&$1.54657$&$0.0319059$&$0.022589$&$0.00044202$&$0.000312945$\\
\hline
$(5,5)$&$0.183437$&$0.129871$&$0.000930408$&$0.000658718$&$0.0000373219$&$0.0000264235$\\
\hline
$(6,4)$&$1.80936$&$1.281$&$0.0294761$&$0.0208687$&$0.000584608$&$0.000413895$\\
\hline
\end{tabular}
\end{center}
\end{minipage}}
\label{tab:multicol}
\end{table}

\subsection{Proofs}
{\bf Proof of Theorem \ref{aoop}}. Differentiating (\ref{wtpgf}) w.r.t. $z$, for $\left|z-a(\bar{p}) z^k\right|<1$, we get
\begin{align*}
M_{\acute{T}}^{\prime}(z) &= n M_{\acute{T}}(z) \frac{(1-k a(\bar{p})z^{k-1})}{(1-z+a(\bar{p})z^k)}=n M_{\acute{T}}(z) (1-k a(\bar{p})z^{k-1}) \sum_{m=0}^{\infty}\left(\sum_{l=0}^{\floor{m/k}}(-1)^l{m-l(k-1) \choose l} a(\bar{p})^l\right)z^m\\
&= n \left(\sum_{m=0}^{\infty} \check{p}_m z^m\right) \left\{\sum_{m=0}^{\infty}b_{m,\bar{p}}z^m-k a(\bar{p}) \sum_{m=k-1}^{\infty}b_{m-k+1,\bar{p}}z^m\right\}\\
&= n  \left\{\sum_{m=0}^{\infty}\left(\sum_{l=0}^{m} \check{p}_l b_{m-l,\bar{p}} \right)z^m-k a(\bar{p}) \sum_{m=k-1}^{\infty}\left(\sum_{l=0}^{m-k+1} \check{p}_l b_{m-k+1-l,\bar{p}} \right)z^m\right\},
\end{align*}
where $\check{p}_m = {\bf P}(\acute{T}=m)$ and $b_{m,\bar{p}}$ defined as in (\ref{nots}). Multiplying by $(1-q z)$ and collecting the coefficients of $z^m$, we get the recurrence relation
\begin{align*}
q \left(\frac{n}{q}+m\right) \check{p}_{m}-(m+1) \check{p}_{m+1} &+n \left[\sum_{l=0}^{m-1}\check{p}_l\left( b_{m-l,\bar{p}}- q b_{m-l-1,\bar{p}}\right)\right.\\
&\left.~~~~~~~~~~~~~~~~~~~-k a(\bar{p})\sum_{l=0}^{m-k+1}\check{p}_l b_{m-k+1-l,\bar{p}}+ q k a(\bar{p})\sum_{l=0}^{m-k}\check{p}_l b_{m-k-l,\bar{p}}\right]=0,
\end{align*}
where $q=1-p \in (0,1)$ defined as in (\ref{Mmean2}). Let $g \in {\cal G}_{\acute{T}}$, defined in (\ref{gx}), then
\begin{align*}
\sum_{m=0}^{\infty} g(m+1)&\left\{q \left(\frac{n}{q}+m\right) \check{p}_{m}-(m+1) \check{p}_{m+1} +n\sum_{l=0}^{m-1}\check{p}_l\left( b_{m-l,\bar{p}}- q b_{m-l-1,\bar{p}}\right)\right.\\
&~~~~~~~~~~~~~~~~~~~~~~~\left.-n k a(\bar{p})\sum_{l=0}^{m-k+1}\check{p}_l b_{m-k+1-l,\bar{p}}+ n q k a(\bar{p})\sum_{l=0}^{m-k}\check{p}_l b_{m-k-l,\bar{p}}\right\}=0.
\end{align*}
This leads to the following
\begin{align*}
\sum_{m=0}^{\infty} &\left[q \left(\frac{n}{q}+m\right) g(m+1)-m g(m) +n\sum_{l=1}^{\infty}g(m+l+1)\left( b_{l,\bar{p}}- q b_{l-1,\bar{p}}\right)\right.\\
&~~~~~~~~~~~~~\left.-n k a(\bar{p})\sum_{l=k-1}^{\infty}g(m+l+1) b_{l-k+1,\bar{p}}+ n q k a(\bar{p})\sum_{l=k}^{\infty} g(m+l+1) b_{l-k,\bar{p}}\right]\check{p}_m=0.
\end{align*}
Hence, Stein operator of $\acute{T}$ is given by  
\begin{align}
{\cal A}_{\acute{T}} g(m) &=q \left(\frac{n}{q}+m\right) g(m+1)-m g(m) +n\sum_{l=1}^{\infty}g(m+l+1)\left( b_{l,\bar{p}}- q b_{l-1,\bar{p}}\right)\nonumber\\
&~~~~~~~~~~~~~~~-n k a(\bar{p})\sum_{l=k-1}^{\infty}g(m+l+1) b_{l-k+1,\bar{p}}+ n q k a(\bar{p})\sum_{l=k}^{\infty} g(m+l+1) b_{l-k,\bar{p}}\label{opso}\\
&= \bar{\cal A}_Z g(m)+ {\cal U}_{\acute{T}}g(m),\nonumber
\end{align}
where $\bar{\cal A}_Z$ denote Stein operator for NB$(\frac{n}{q},p)$. This is a Stein operator for $\acute{T}$ in perturbation of NB operator. Using (\ref{delta}) in perturbed operator ${\cal U}_{\acute{T}}$, we get 
\begin{align*}
{\cal U}_{\acute{T}}  g(m) \hspace{-0.1cm}=& n \hspace{-0.1cm}\left[\sum_{l=1}^{\infty}b_{l,\bar{p}}- q \sum_{l=1}^{\infty}b_{l-1,\bar{p}}-k a(\bar{p}) \sum_{l=k-1}^{\infty}b_{l-k+1,\bar{p}} + q k a(\bar{p})\sum_{l=k}^{\infty}b_{l-k,\bar{p}}\right]\hspace{-0.1cm}g(m\hspace{-0.1cm}+\hspace{-0.1cm}1)+n\sum_{l=1}^{\infty} \sum_{j=1}^{l} \Delta g(m\hspace{-0.1cm}+\hspace{-0.1cm}j) b_{l,\bar{p}} \\
&- nq \sum_{l=1}^{\infty} \sum_{j=1}^{l} \Delta g(m\hspace{-0.1cm}+\hspace{-0.1cm}j) b_{l-1,\bar{p}} - nk a(\bar{p}) \hspace{-0.2cm}\sum_{l=k-1}^{\infty} \sum_{j=1}^{l} \Delta g(m\hspace{-0.1cm}+\hspace{-0.1cm}j) b_{l-k+1,\bar{p}}+nq k a(\bar{p}) \sum_{l=k}^{\infty} \sum_{j=1}^{l} \Delta g(m\hspace{-0.1cm}+\hspace{-0.1cm}j) b_{l-k,\bar{p}}
\end{align*}
Observe that $\sum_{l=0}^{\infty} b_{l,\bar{p}} = 1/a(\bar{p})$ and $b_{0,\bar{p}} = 1$. Using (\ref{Mmean2}), we obtain the perturbation operator
\begin{align*}
{\cal U}_{\acute{T}}  g(m) &= n\sum_{l=1}^{\infty} \sum_{j=1}^{l} \Delta g(m+j) b_{l,\bar{p}} - nq \sum_{l=1}^{\infty} \sum_{j=1}^{l} \Delta g(m+j) b_{l-1,\bar{p}}\\
& ~~~~~~~~~~~~~~~~~~~~~~~~~~~~~~~- nk a(\bar{p}) \sum_{l=k-1}^{\infty} \sum_{j=1}^{l} \Delta g(m+j) b_{l-k+1,\bar{p}}+nq k a(\bar{p}) \sum_{l=k}^{\infty} \sum_{j=1}^{l} \Delta g(m+j) b_{l-k,\bar{p}}
\end{align*}
Hence, for $g \in {\cal G}_Z \cap {\cal G}_{\acute{T}} $, taking expectation w.r.t $\acute{T}$ and using (\ref{bound}), we get required result.\qed

\noindent
{\bf Proof of Theorem \ref{atop}}. Next, for $g \in {\cal G}_{\acute{T}} $, we Introduce a new parameter $\alpha > 0$ in (\ref{opso}) as
\begin{align*}
{\cal A}_{\acute{T}}  g(m) &=q \left(\alpha+m\right) g(m+1)-m g(m)+(n-\alpha q)g(m+1) +n\sum_{l=1}^{\infty}g(m+l+1)\left( b_{l,\bar{p}}- q b_{l-1,\bar{p}}\right)\\
&~~~~~~~~~~~~~~~~~~~~~~~~~~~~~~~~~~~~~-n k a(\bar{p})\sum_{l=k-1}^{\infty}g(m+l+1) b_{l-k+1,\bar{p}}+ n q k a(\bar{p})\sum_{l=k}^{\infty} g(m+l+1) b_{l-k,\bar{p}}\\
&={\cal A}_Z g(m) + \hat{\cal U}_{\acute{T}} g(m)\label{tpo}
\end{align*}
where $\alpha$ and $q=1-p$ as defined in (\ref{alphap}). This is Stein operator for $\acute{T}$, which is a perturbation of NB$(\alpha, p)$. 
Putting (\ref{delta}) and (\ref{ne}) in perturbed operator $\hat{\cal U}_{\acute{T}}$ then using(\ref{mvar}), we get
\begin{align*}
\hat{\cal U}_{\acute{T}}  g(m)&=n \left[\sum_{l=1}^{\infty} \sum_{j=1}^{l} \sum_{u=1}^{j-1}\Delta^2 g(m+u) b_{l,\bar{p}}- q \sum_{l=1}^{\infty} \sum_{j=1}^{l}\sum_{u=1}^{j-1} \Delta^2 g(m+u) b_{l-1,\bar{p}}\right.\\
&~~~~~~~~~~~~~~~~~~~~~\left.-k a(\bar{p}) \sum_{l=k-1}^{\infty} \sum_{j=1}^{l}\sum_{u=1}^{j-1} \Delta^2 g(m+u) b_{l-k+1,\bar{p}} + k q a(\bar{p})\sum_{l=k}^{\infty} \sum_{j=1}^{l} \sum_{u=1}^{j-1}\Delta^2 g(m+u) b_{l-k,\bar{p}}\right].
\end{align*}
Taking the expectation w.r.t. $\acute{T}$, we have
\begin{align*}
{\mathbb E}[\hat{\cal U}_{\acute{T}}  g(\acute{T})]&= n\sum_{m=0}^{\infty} \left[\sum_{l=1}^{\infty} \sum_{j=1}^{l} \sum_{u=1}^{j-1}\Delta^2 g(m+u) b_{l,\bar{p}}- q \sum_{l=1}^{\infty} \sum_{j=1}^{l}\sum_{u=1}^{j-1} \Delta^2 g(m+u) b_{l-1,\bar{p}}\right.\\
&~~~~~~~~~~~~~~~\left.-k a(\bar{p}) \sum_{l=k-1}^{\infty} \sum_{j=1}^{l}\sum_{u=1}^{j-1} \Delta^2 g(m+u) b_{l-k+1,\bar{p}} + k q a(\bar{p})\sum_{l=k}^{\infty} \sum_{j=1}^{l} \sum_{u=1}^{j-1}\Delta^2 g(m+u) b_{l-k,\bar{p}}\right]\check{p}_m\\
&= n \sum_{m=0}^{\infty}\left[\sum_{l=1}^{\infty} \sum_{j=1}^{l} \sum_{u=1}^{j-1}\Delta g(m+u)   b_{l,\bar{p}}- q \sum_{l=1}^{\infty} \sum_{j=1}^{l}\sum_{u=1}^{j-1} \Delta g(m+u)  b_{l-1,\bar{p}}\right.\\
&~~\left.-k a(\bar{p}) \sum_{l=k-1}^{\infty} \sum_{j=1}^{l}\sum_{u=1}^{j-1} \Delta g(m+u) b_{l-k+1,\bar{p}} + k q a(\bar{p})\sum_{l=k}^{\infty} \sum_{j=1}^{l} \sum_{u=1}^{j-1}\Delta g(m+u)  b_{l-k,\bar{p}}\right](\check{p}_{m-1}-\check{p}_{m})
\end{align*}
Now, $g \in {\cal G}_Z \cap {\cal G}_{\acute{T}} $, taking supremum and using (\ref{bound}), we get required result.\qed

\singlespacing

\end{document}